\documentclass[a4paper, twoside]{article}

\usepackage[mathscr]{eucal}
\usepackage{amsmath,amsfonts,amssymb,amsthm}
\usepackage{mathrsfs}
\usepackage[dvips]{graphicx}
\usepackage{fancyhdr}

\newcommand{\mathsym}[1]{{}}

\newtheorem{theorem}{\mdseries\scshape Theorem}[section]

\newtheorem{proposition}{\mdseries\scshape Proposition}[section]
\newtheorem{lemma}{\mdseries\scshape Lemma}[section]

\newtheorem{conjecture}{\mdseries\scshape Conjecture}[section]

\renewcommand{\footrulewidth}{0pt}

\begin{document}
\title{\large\bfseries ON THE ZEROS OF THE EISENSTEIN SERIES \\ FOR $\Gamma_0^* (5)$ AND $\Gamma_0^* (7)$}
\author{Junichi Shigezumi}
\date{}

\maketitle

\noindent
{\bfseries Abstract.}
We locate almost all the zeros of the Eisenstein series associated with the Fricke groups of level $5$ and $7$ in their fundamental domains by applying and extending the method of F. K. C. Rankin and H. P. F. Swinnerton-Dyer (1970). We also use the arguments of some terms of the Eisenstein series in order to improve existing error bounds.

\thispagestyle{fancy}
\lhead{}\chead{}\rhead{}
\lfoot{{\small {\bfseries 2000 Mathematics Subject Classification:} Primary 11F11.}}
\renewcommand{\footrulewidth}{0.4pt}

\section{Introduction}

F. K. C. Rankin and H. P. F. Swinnerton-Dyer considered the problem of locating the zeros of the Eisenstein series $E_k(z)$ in the standard fundamental domain $\mathbb{F}$ \cite{RSD}. They proved that all of the zeros of $E_k(z)$ in $\mathbb{F}$ lie on the unit circle. They also stated towards the end of their study that ``This method can equally well be applied to Eisenstein series associated with subgroups of the modular group.'' However, it seems unclear how widely this claim holds. 

Subsequently, T. Miezaki, H. Nozaki, and the present author considered the same problem for the Fricke group $\Gamma_0^{*}(p)$ (See \cite{K}, \cite{Q}), and proved that all of the zeros of the Eisenstein series $E_{k, p}^{*}(z)$ in a certain fundamental domain lie on a circle whose radius is equal to $1 / \sqrt{p}$,  $p = 2, 3$ \cite{MNS}.

The Fricke group $\Gamma_0^{*}(p)$ is not a subgroup of $\text{SL}_2(\mathbb{Z})$, but it is commensurable with $\text{SL}_2(\mathbb{Z})$. For a fixed prime $p$, we define $\Gamma_0^{*}(p) := \Gamma_0(p) \cup \Gamma_0(p) \: W_p$, where $\Gamma_0(p)$ is a congruence subgroup of $\text{SL}_2(\mathbb{Z})$.

Let $k \geqslant 4$ be an even integer. For $z \in \mathbb{H} := \{z \in \mathbb{C} \: ; \: Im(z)>0 \}$, let
\begin{equation}
E_{k, p}^{*}(z) := \frac{1}{p^{k / 2}+1} \left(p^{k/2} E_k(p z) + E_k(z) \right) \label{def:e*}
\end{equation}
be the Eisenstein series associated with $\Gamma_0^{*}(p)$. ({\it cf.} \cite{SG})

\renewcommand{\thefootnote}{\fnsymbol{footnote}}

Henceforth, we assume that $p = 5$ or $7$. The region\footnote[2]{In the paper published on Kyushu J. Math. ({\bfseries 61}(2007), 527--549), there is a mistake on the definition of $\mathbb{F}^{*}(p)$. The definition in this paper is correct. We thank Prof. Rainer Schulze-Pillot for pointing it out.}
\begin{align}
\mathbb{F}^{*}(p) &:= \left\{|z| \geqslant 1 / \sqrt{p}, \: |z + 1/2| \geqslant 1 / (2 \sqrt{p}), \: - 1 / 2 \leqslant Re(z) \leqslant 0\right\} \notag\\
 &\qquad \qquad \bigcup \left\{|z| > 1 / \sqrt{p}, \: |z - 1/2| > 1 / (2 \sqrt{p}), \: 0 \leqslant Re(z) < 1 / 2 \right\} \label{fd-0sp}
\end{align}
is a fundamental domain for $\Gamma_0^{*}(p)$. ({\it cf.} \cite{SH}, \cite{SE}) Define $A_p^{*} := \mathbb{F}^{*}(p) \cap \{z \in \mathbb{C} \: ; \: |z| = 1 / \sqrt{p} \; \text{or} \; |z \pm 1/2| = 1 / (2 \sqrt{p})\}$.

In the present paper, we will apply the method of F. K. C. Rankin and H. P. F. Swinnerton-Dyer (RSD Method) to the Eisenstein series associated with $\Gamma_0^{*}(5)$ and $\Gamma_0^{*}(7)$. We have the following conjectures:

\begin{conjecture}
Let $k \geqslant 4$ be an even integer. Then all of the zeros of $E_{k, 5}^{*}(z)$ in $\mathbb{F}^{*}(5)$ lie on the arc $A_5^{*}$. \label{conj-g0s5}
\end{conjecture}

\begin{conjecture}
Let $k \geqslant 4$ be an even integer. Then all of the zeros of $E_{k, 7}^{*}(z)$ in $\mathbb{F}^{*}(7)$ lie on the arc $A_7^{*}$. \label{conj-g0s7}
\end{conjecture}

First, we prove that all but at most $2$ zeros of $E_{k, p}^{*}(z)$ in $\mathbb{F}^{*}(p)$ lie on the arc $A_p^{*}$ (See Subsection \ref{subsec-g0s5-ab2} and \ref{subsec-g0s7-ab2}). Second, if $(24 / (p + 1)) \mid k$, we prove that all of the zeros of $E_{k, p}^{*}(z)$ in $\mathbb{F}^{*}(p)$ lie on $A_p^{*}$ (See Subsection \ref{subsec-g0s5-40} and \ref{subsec-g0s7-60}).

We can then prove that if $(24 / (p + 1)) \nmid k$, all but one of the zeros of $E_{k, p}^{*}(z)$ in $\mathbb{F}^{*}(p)$ lie on $A_p^{*}$. Furthermore, let $\alpha_5 \in [0, \pi]$ (resp. $\alpha_7 \in [0, \pi]$) be the angle which satisfies $\tan\alpha_5 = 2$ (resp. $\tan\alpha_7 = 5 / \sqrt{3}$), and let $\alpha_{p, k} \in [0, \pi]$ be the angle which satisfies $\alpha_{p, k} \equiv k (\pi / 2 + \alpha_p) / 2 \pmod{\pi}$. Then, since $\alpha_p$ is an irrational multiple of $\pi$, $\alpha_{p, k}$ appear in the interval $[0, \pi]$ uniformly for all even integers $k \geqslant 4$. In Subsection \ref{subsec-g0s5-41}, we prove that all of the zeros of $E_{k, 5}^{*}(z)$ in $\mathbb{F}^{*}(5)$ are on $A_5^{*}$ if $\alpha_{5, k} < (116/180) \pi$ or $(117/180) \pi < \alpha_{5, k}$. That is, we prove about $179/180$ of Conjecture \ref{conj-g0s5}. Similarly, in Subsection \ref{subsec-g0s7-61} and \ref{subsec-g0s7-62}, we prove that all of the zeros of $E_{k, 7}^{*}(z)$ in $\mathbb{F}^{*}(7)$ are on $A_7^{*}$ if ``$\alpha_{7, k} < (127.68/180) \pi$ or $(128.68/180) \pi < \alpha_{7, k}$ for $k \equiv 2 \pmod{6}$'' or ``$\alpha_{7, k} < (108.5/180) \pi$ or $(109.5/180) \pi < \alpha_{7, k}$ for $k \equiv 4 \pmod{6}$''. Thus  we can also prove about $179/180$ of Conjecture \ref{conj-g0s7}.

In \cite{RSD}, we considered a bound for the error terms $R_1$ (See (\ref{eq-fkt})) in terms only of their absolute values. However, in the present paper, we also use the arguments of some terms in the series.  We can then approach the exact value of the Eisenstein series.

A more detailed account of the material in the present study may be found in \cite{SJ}.

\section{General Theory}

\lhead[\thepage]{\quad}\rhead[\quad]{\thepage}
\chead[{\upshape J. Shigezumi}]{\upshape Eisenstein series for $\Gamma_0^* (5)$ and $\Gamma_0^* (7)$}
\lfoot{}
\renewcommand{\footrulewidth}{0pt}

\subsection{Preliminaries}
Let $v_p(f)$ be the order of a modular function $f$ at a point $p$.

\subsubsection{$\Gamma_0^{*}(5)$}
We define
\begin{align*}
A_{5, 1}^{*} &:= \{ z \: ; \: |z| = 1 / \sqrt{5}, \: \pi / 2 < Arg(z) < \pi / 2 + \alpha_5 \},\\
A_{5, 2}^{*} &:= \{ z \: ; \: |z + 1 / 2| = 1 / (2 \sqrt{5}), \: \alpha_5 < Arg(z) < \pi / 2 \}.
\end{align*}
Then, $A_5^{*} = A_{5, 1}^{*} \cup A_{5, 2}^{*} \cup \{ i / \sqrt{5}, \: \rho_{5, 1}, \: \rho_{5, 2} \}$, where $\rho_{5, 1} := - 1 / 2 + i / \left(2 \sqrt{5}\right)$ and $\rho_{5, 2} := - 2 / 5 + i / 5$.

Let $f$ be a modular form for $\Gamma_0^{*}(5)$ of weight $k$, and let $k$ be an even integer such that $k \equiv 2 \pmod{4}$, then
\begin{equation*}
f(i / \sqrt{5}) = f(W_5 \: i / \sqrt{5}) = i^k f(i / \sqrt{5}) = - f(i / \sqrt{5}).
\end{equation*}
Thus, we have $f(i / \sqrt{5}) = 0$. Similarly, we have $f(\rho_{5, 1}) = f(\rho_{5, 2}) = 0$. Thus, we have $v_{i / \sqrt{5}} (f) \geqslant 1$, $v_{\rho_{5, 1}} (f) \geqslant 1$, and $v_{\rho_{5, 2}} (f) \geqslant 1$.

On the other hand, if $k \equiv 0 \pmod{4}$, then we have $v_{i / \sqrt{5}} (E_{k, 5}^{*}) = v_{\rho_{5, 1}} (E_{k, 5}^{*}) = v_{\rho_{5, 2}} (E_{k, 5}^{*}) = 0$.

\subsubsection{$\Gamma_0^{*}(7)$}
We define
\begin{align*}
A_{7, 1}^{*} &:= \{ z \: ; \: |z| = 1 / \sqrt{7}, \: \pi / 2 < Arg(z) < \pi / 2 + \alpha_7 \},\\
A_{7, 2}^{*} &:= \{ z \: ; \: |z + 1 / 2| = 1 / (2 \sqrt{7}), \: \alpha_7 - \pi / 6 < Arg(z) < \pi / 2 \}.
\end{align*}
Then, $A_7^{*} = A_{7, 1}^{*} \cup A_{7, 2}^{*} \cup \{ i / \sqrt{7}, \: \rho_{7, 1}, \: \rho_{7, 2} \}$, where $\rho_{7, 1} := - 1 / 2 + i / \left(2 \sqrt{7}\right)$ and $\rho_{7, 2} := - 5 / 14 + \sqrt{3} i / 14$. 

Let $f$ be a modular form for $\Gamma_0^{*}(7)$ of weight $k$. If $k \equiv 2 \pmod{4}$, then we have $v_{i / \sqrt{7}} (f) \geqslant 1$ and $v_{\rho_{7, 1}} (f) \geqslant 1$. On the other hand, if $k \equiv 0 \pmod{4}$, then we have $v_{i / \sqrt{7}} (E_{k, 7}^{*}) = v_{\rho_{7, 1}} (E_{k, 7}^{*}) = 0$.

Similarly, if $k \not\equiv 0 \pmod{6}$, then we have $v_{\rho_{7, 2}} (f) \geqslant 1$, while if $k \equiv 0 \pmod{6}$, then we have $v_{\rho_{7, 2}} (E_{k, 7}^{*}) = 0$.

\subsection{Valence Formula}
In order to determine the location of zeros of $E_{k,p}^{*}(z)$ in $\mathbb{F}^{*}(p)$, we require the valence formula for $\Gamma_0^{*}(p)$.

\begin{proposition}\label{prop-vf-g0s5}
Let $f$ be a modular function of weight $k$ for $\Gamma_0^{*}(5)$, which is not identically zero. We have
\begin{equation}
v_{\infty}(f) + \frac{1}{2} v_{i / \sqrt{5}}(f) + \frac{1}{2} v_{\rho_{5,1}} (f) + \frac{1}{2} v_{\rho_{5,2}} (f) + \sum_{\begin{subarray}{c} p \in \Gamma_0^{*}(5) \setminus \mathbb{H} \\ p \ne i / \sqrt{5}, \rho_{5, 1}, \rho_{5, 2}\end{subarray}} v_p(f) = \frac{k}{4}.
\end{equation}
\end{proposition}

\begin{proposition}\label{prop-vf-g0s7}
Let $f$ be a modular function of weight $k$ for $\Gamma_0^{*}(7)$, which is not identically zero. We have
\begin{equation}
v_{\infty}(f) + \frac{1}{2} v_{i / \sqrt{7}}(f) + \frac{1}{2} v_{\rho_{7,1}} (f) + \frac{1}{3} v_{\rho_{7,2}} (f) + \sum_{\begin{subarray}{c} p \in \Gamma_0^{*}(7) \setminus \mathbb{H} \\ p \ne i / \sqrt{7}, \rho_{7, 1}, \rho_{7, 2}\end{subarray}} v_p(f) = \frac{k}{3}.
\end{equation}
\end{proposition}

The proofs of the above propositions are very similar to that for the valence formula for $\text{SL}_2(\mathbb{Z})$ ({\it cf.} \cite{SE}).

\subsection{Some Eisenstein series of low weights}

By means of a straightforward calculation, we have the following propositions:

\begin{proposition}
The location of the zeros of the Eisenstein series $E_{k, 5}^{*}$ in $\mathbb{F}^{*}(5)$ for $4 \geqslant k \geqslant 10$ are given by the following table:
\begin{center}
\begin{tabular}{ccccccc}
$k$ & $v_{\infty}$ & $v_{i / \sqrt{5}}$ & $v_{\rho_{5, 1}}$ & $v_{\rho_{5, 2}}$ & $V_{5, 1}^{*}$ & $V_{5, 2}^{*}$\\
\hline
$4$ & $0$ & $0$ & $0$ & $0$ & $1$ & $0$\\
$6$ & $0$ & $1$ & $1$ & $1$ & $0$ & $0$\\
$8$ & $0$ & $0$ & $0$ & $0$ & $1$ & $1$\\
$10$ & $0$ & $1$ & $1$ & $1$ & $1$ & $0$\\
\hline
\end{tabular}
\end{center}
where $V_{5, n}^{*}$ denotes the number of simple zeros of the Eisenstein series $E_{k, 5}^{*}$ on the arc $A_{5, n}^{*}$ for $n = 1, 2$.
\end{proposition}

\begin{proposition}
The location of the zeros of the Eisenstein series $E_{k, 7}^{*}$ in $\mathbb{F}^{*}(7)$ for $k = 4, 6,$ and $12$ are given by the following table: 
\begin{center}
\begin{tabular}{cccccc}
$k$ & $v_{\infty}$ & $v_{i / \sqrt{7}}$ & $v_{\rho_{7, 1}}$ & $v_{\rho_{7, 2}}$ & $V_7^{*}$\\
\hline
$4$ & $0$ & $0$ & $0$ & $1$ & $1$\\
$6$ & $0$ & $1$ & $1$ & $0$ & $1$\\
$12$ & $0$ & $0$ & $0$ & $0$ & $4$\\
\hline
\end{tabular}
\end{center}
where $V_7^{*}$ denotes the number of simple zeros of the Eisenstein series $E_{k, 7}^{*}$ on $A_{7, 1}^{*} \cup A_{7, 2}^{*}$.
\end{proposition}

\subsection{The space of modular forms}

Let $M_{k, p}^{*}$ be the space of modular forms for $\Gamma_0^{*} (p)$ of weight $k$, and let $M_{k, p}^{* 0}$ be the space of cusp forms for $\Gamma_0^{*} (p)$ of weight $k$. Upon considering the map $M_{k, p}^{*} \ni f \mapsto f(\infty) \in \mathbb{C}$, it is clear that the kernel of this map is given by $M_{k, p}^{* 0}$. So $\dim(M_{k, p}^{*} / M_{k, p}^{* 0}) \leqslant 1$, and $M_{k, p}^{*} = \mathbb{C} E_{k, p}^{*} \oplus M_{k, p}^{* 0}$. Let $\eta(z)$ be the {\it Dedekind's $\eta$-function}.

\subsubsection{$\Gamma_0^{*}(5)$}
Note that $\Delta_5 = \eta^4(z) \eta^4(5 z)$ is a cusp form for $\Gamma_0^{*} (5)$ of weight $4$. We have the following theorem:

\begin{theorem}Let $k$ be an even integer.\\
\quad$(1)$ For $k < 0$ and $k = 2$, $M_{k, 5}^{*} = 0$.\\
\quad$(2)$ For $k = 0$ and $6$, we have $M_{k, 5}^{* 0} = 0$, and $\dim(M_{k, 5}^{*}) = 1$ with a base $E_{k, 5}^{*}$.\\
\quad$(3)$ $M_{k, 5}^{* 0} = \Delta_5 M_{k - 4, 5}^{*}$.\label{th-mod_sp_5}
\end{theorem}
The proof of this theorem is very similar to that for $\text{SL}_2(\mathbb{Z})$. Furthermore, we have
\begin{align*}
M_{4 n, 5}^{*} &= \mathbb{C} (E_{4, 5}^{*})^n \oplus \mathbb{C} (E_{4, 5}^{*})^{n - 1} \Delta_5 \oplus \cdots \oplus \mathbb{C} \Delta_5^n,\\
M_{4 n + 6, 5}^{*} &= E_{6, 5}^{*} ((E_{4, 5}^{*})^n \oplus \mathbb{C} (E_{4, 5}^{*})^{n - 1} \Delta_5 \oplus \cdots \oplus \mathbb{C} \Delta_5^n).
\end{align*}
Thus, we have the following proposition:
\begin{proposition}\label{prop-bd_ord_5}
Let $k \geqslant 4$ be an even integer. For every $f \in M_{k, 5}^{*}$, we have
\begin{equation}
\begin{split}
v_{i / \sqrt{5}}(f) \geqslant s_k, \quad v_{\rho_{5, 1}}(f) \geqslant s_k, \quad v_{\rho_{5, 2}}(f) \geqslant s_k\\
(s_k=0, 1 \; \text{such that} \; 2 s_k \equiv k \pmod{4}).
\end{split}
\end{equation}
\end{proposition}

\subsubsection{$\Gamma_0^{*}(7)$}
We define $\Delta_7 := \eta^6(z) \eta^6(7 z)$ and ${E_{2, 7}}'(z) := (7 E_2(7 z) - E_2(z)) / 6$. We then have the following theorem:
\begin{theorem}Let $k$ be an even integer, and we define $\Delta_{7, 4} := (5 / 16) (({E_{2, 7}}')^2 - E_{4, 7}^{*})$ and\\ \quad $\Delta_{7, 10}^{0} := (559 / 690) ((41065 / 137592) (E_{4, 7}^{*} E_{6, 7}^{*} - E_{10, 7}^{*}) - E_{6, 7}^{*} \Delta_{7, 4})$\\
\quad$(1)$ $M_{k, 7}^{* 0} = M_{12, 7}^{* 0} M_{k - 12, 7}^{*}$.\\
\quad$(2)$ For $k < 0$ and $k = 2$, $M_{k, 7}^{*} = 0$. We have $M_{0, 7}^{*} = \mathbb{C}$.\\
\quad$(3)$ We have $M_{4, 7}^{* 0} = \mathbb{C} \Delta_{7, 4}$, $M_{6, 7}^{* 0} = \mathbb{C} \Delta_{7, 10}^0 / \Delta_{7, 4}$,\\
\qquad $M_{8, 7}^{* 0} = \mathbb{C} (\Delta_{7, 4})^2 \oplus \mathbb{C} E_{4, 7}^{*} \Delta_{7, 4}$, $M_{10, 7}^{* 0} = \mathbb{C} \Delta_{7, 10}^0 \oplus \mathbb{C} E_{6, 7}^{*} \Delta_{7, 4}$,\\
\qquad $M_{12, 7}^{* 0} = \mathbb{C} (\Delta_7)^2 \oplus \mathbb{C} (\Delta_{7, 4})^3 \oplus \mathbb{C} E_{4, 7}^{*} (\Delta_{7, 4})^2 \oplus \mathbb{C} (E_{4, 7}^{*})^2 \Delta_{7, 4}$, and\\
\qquad $M_{14, 7}^{* 0} = \mathbb{C} \Delta_{7, 4} \Delta_{7, 10}^0 \oplus \mathbb{C} E_{6, 7}^{*} (\Delta_{7, 4})^2 \oplus \mathbb{C} E_{4, 7}^{*} E_{6, 7}^{*} \Delta_{7, 4}$.\label{th-mod_sp_7}
\end{theorem}

 We thus have the following table:

\begin{allowdisplaybreaks}
\begin{center}
\begin{tabular}{rcccccc}
$k$ & $f$ & $v_{\infty}$ & $v_{i / \sqrt{7}}$ & $v_{\rho_{7, 1}}$ & $v_{\rho_{7, 2}}$ & $V_7^{*}$\\
\hline
$4$ & $({E_{2, 7}}')^2$ & $0$ & $0$ & $0$ & $4$ & $0$\\
 & $\Delta_{7, 4}$ & $1$ & $0$ & $0$ & $1$ & $0$\\
$10$ & $\Delta_{7, 10}^0$ & $2$ & $1$ & $1$ & $1$ & $0$\\
$12$ & $(\Delta_7)^2$ & $4$ & $0$ & $0$ & $0$ & $0$\\
\hline
\end{tabular}
\end{center}
\end{allowdisplaybreaks}

Furthermore, we have
\begin{align*}
M_{k, 7}^{*}
 &= E_{k - 12 n, 7}^{*} \left\{\mathbb{C} (E_{4, 7}^{*})^{3 n} \oplus (E_{4, 7}^{*})^{3 (n - 1)} M_{12, 7}^{* 0} \oplus (E_{4, 7}^{*})^{3 (n - 2)} (M_{12, 7}^{* 0})^2 \oplus \cdots \oplus (M_{12, 7}^{* 0})^n \right\}\\
 &\qquad \oplus M_{k - 12 n, 7}^{* 0} (M_{12, 7}^{* 0})^n.
\end{align*}
Thus, we have the following proposition:
\begin{proposition}
Let $k \geqslant 4$ be an even integer. For every $f \in M_{k, 7}^{*}$, we have
\begin{equation}\label{prop-bd_ord_7}
\begin{split}
v_{i / \sqrt{7}}(f) \geqslant s_k, \quad v_{\rho_{7, 1}}(f) \geqslant s_k \quad
 &(s_k=0, 1 \; \text{such that} \; 2 s_k \equiv k \pmod{4}),\\
v_{\rho_{7, 2}}(f) \geqslant t_k \quad
 &(s_k=0, 1, 2 \; \text{such that} \; - 2 t_k \equiv k \pmod{6}).
\end{split}
\end{equation}
\end{proposition}

\section{The method of Rankin and Swinnerton-Dyer (RSD Method)}

\subsection{RSD Method}
Let $k \geqslant 4$ be an even integer. For $z \in \mathbb{H}$, we have
\begin{equation}
E_k(z) = \frac{1}{2} \sum_{(c,d)=1}(c z + d)^{-k}. \label{def:e}
\end{equation}
Moreover, we have $\mathbb{F} = \left\{|z| \geqslant 1, \: - 1 / 2 \leqslant Re(z) \leqslant 0\right\} \cup \left\{|z| > 1, \: 0 \leqslant Re(z) < 1 / 2 \right\}$.

At the beginning of their proof in \cite{RSD},  Rankin and Swinnerton-Dyer considered the following function:
\begin{equation}
F_k(\theta) := e^{i k \theta / 2} E_k\left(e^{i \theta}\right), \label{def:f}
\end{equation}
which is real for all $\theta \in [0, \pi]$. Considering the four terms with $c^2 + d^2 = 1$, they proved that
\begin{equation}
F_k(\theta) = 2 \cos(k \theta / 2) + R_1, \label{eq-fkt}
\end{equation}
where $R_1$ denotes the remaining terms of the series. Moreover they showed that $|R_1| < 2$ for all $k \geqslant 12$.  If $\cos (k \theta / 2)$ is $+1$ or $-1$, then $F_k(2 m \pi / k)$ is positive or negative, respectively, and we can show the existence of the zeros. In addition, we can prove that all of the zeros of $E_k(z)$ in $\mathbb{F}$ lie on the unit circle using the {\it Valence Formula} and the theory on the space of modular forms for $\text{SL}_2(\mathbb{Z})$.

\subsection{The function: $F_{k, p, n}^{*}$}
We expect all of the zeros of the Eisenstein series $E_{k, p}^{*}(z)$ in $\mathbb{F}^{*}(p)$ to lie on the arcs $e^{i \theta} / \sqrt{p}$ and $e^{i \theta} / (2 \sqrt{p}) - 1 / 2$, which form the boundary of the fundamental domain defined by the equation (\ref{fd-0sp}).

We define
\begin{gather}
F_{k, p, 1}^{*}(\theta) := e^{i k \theta / 2} E_{k,p}^{*}\left(e^{i \theta} / \sqrt{p}\right), \label{def:f*1}\\
F_{k, p, 2}^{*}(\theta) := e^{i k \theta / 2} E_{k,p}^{*}\left(e^{i \theta} / 2 \sqrt{p} - 1 / 2\right). \label{def:f*p2}
\end{gather}

We can write
\begin{equation}
F_{k, p, 1}^{*}(\theta)
 = \frac{1}{2} \sum_{\begin{subarray}{c} (c,d)=1\\ p \nmid c\end{subarray}}(c e^{i \theta / 2} + \sqrt{p} d e^{-i \theta / 2})^{-k}
 + \frac{1}{2} \sum_{\begin{subarray}{c} (c,d)=1\\ p \nmid c\end{subarray}}(c e^{-i \theta / 2} + \sqrt{p} d e^{i \theta / 2})^{-k},
\end{equation}
\begin{equation}
\begin{split}
F_{k, p, 2}^{*}(\theta)
 &= \frac{1}{2} \sum_{\begin{subarray}{c} (c, d) = 1\\ p \nmid c \; 2 \mid c d \end{subarray}}\left( c e^{i \theta / 2} + d \sqrt{p} e^{- i \theta / 2} \right)^{-k}
+ \frac{1}{2} \sum_{\begin{subarray}{c} (c, d) = 1\\ p \nmid c \; 2 \mid c d \end{subarray}}\left( c e^{- i \theta / 2} + d \sqrt{p} e^{i \theta / 2} \right)^{-k}\\
&+ \frac{2^k}{2} \sum_{\begin{subarray}{c} (c, d) = 1\\ p \nmid c \; 2 \nmid c d \end{subarray}}\left( c e^{i \theta / 2} + d \sqrt{p} e^{- i \theta / 2} \right)^{-k}
+ \frac{2^k}{2} \sum_{\begin{subarray}{c} (c, d) = 1\\ p \nmid c \; 2 \nmid c d \end{subarray}}\left( c e^{- i \theta / 2} + d \sqrt{p} e^{i \theta / 2} \right)^{-k}.
\end{split}
\end{equation}

Hence we can use these expressions as definitions. Note that $(c e^{i \theta / 2} + \sqrt{p} d e^{-i \theta / 2})^{-k}$ and $(c e^{-i \theta / 2} + \sqrt{p} d e^{i \theta / 2})^{-k}$ are conjugates of each other for any pair $(c, \: d)$. Thus, we have the following proposition:

\begin{proposition}
$F_{k, p, 1}^{*}(\theta)$ is real, for all $\theta \in [0, \pi]$. \label{prop-f*1}
\end{proposition}
\begin{proposition}
$F_{k, p, 2}^{*}(\theta)$ is real, for all $\theta \in [0, \pi]$. \label{prop-f*2}
\end{proposition}

Now, we define
\begin{equation*}
F_{k, 5}^{*} (\theta)=
\begin{cases}
F_{k, 5, 1}^{*}(\theta) & \pi / 2 \leqslant \theta \leqslant \pi / 2 + \alpha_5\\
F_{k, 5, 2}^{*}(\theta - \pi / 2) & \pi / 2 + \alpha_5 \leqslant \theta \leqslant \pi
\end{cases}.
\end{equation*}
Then, $F_{k, 5}^{*}$ is continuous in the interval $[\pi / 2, \pi]$. Note that $F_{k, 5, 1}^{*} (\pi / 2 + \alpha_5) = e^{i (\pi / 2) k / 2} F_{k, 5, 2}^{*} (\alpha_5)$.

Similarly,
\begin{equation*}
F_{k, 7}^{*} (\theta)=
\begin{cases}
F_{k, 7, 1}^{*}(\theta) & \pi / 2 \leqslant \theta \leqslant \pi / 2 + \alpha_7\\
F_{k, 7, 2}^{*}(\theta - 2 \pi / 3) & \pi / 2 + \alpha_7 \leqslant \theta \leqslant 7 \pi / 6
\end{cases}.
\end{equation*}
whereupon  $F_{k, 7}^{*}$ is continuous in the interval $[\pi / 2, 7 \pi / 6]$. Note also that $F_{k, 7, 1}^{*} (\pi / 2 + \alpha_7) = e^{i (2 \pi / 3) k / 2} F_{k, 7, 2}^{*} (\alpha_7 - \pi / 6)$.

\subsection{Application of RSD Method}

We introduce $N := c^2 + d^2$. First, we consider the case $N = 1$.  For this case, we can write
\begin{gather}
F_{k, p, 1}^{*}(\theta) = 2 \cos(k \theta /2) + R_{p, 1}^{*},\\
F_{k, p, 2}^{*}(\theta) = 2 \cos(k \theta /2) + R_{p, 2}^{*},
\end{gather}
where $R_{p, 1}^{*}$ and $R_{p, 2}^{*}$ denote the terms satisfying $N > 1$ of $F_{k, p, 1}^{*}$ and $F_{k, p, 2}^{*}$, respectively.

\subsubsection{For $\Gamma_0^{*}(5)$}

For $R_{5, 1}^{*}$, we will consider the following cases: $N = 2$, $5$, $10$, $13$, $17$, and $N \geqslant 25$. Considering $- 2 / \sqrt{5} \leqslant \cos\theta \leqslant 0$, we have
\begin{align}
|R_{5, 1}^{*}|
 &\leqslant 2
 + 4 \left(\frac{1}{2}\right)^{k/2} + 2 \left(\frac{1}{3}\right)^{k/2}
 + \cdots + 2 \left(\frac{1}{9}\right)^{k}
 + \frac{384 \sqrt{6}}{k - 3} \left(\frac{1}{2}\right)^{k}. \label{r*51bound0}
\end{align}

Similarly, for $R_{5, 2}^{*}$, we will consider the cases: $N = 2$, $5$, $10$, $13$, $17$, $25$, $26$, $29$, and $N \geqslant 34$. For these cases we have
\begin{equation}
|R_{5, 2}^{*}|
 \leqslant 2
 + 2 \left(\frac{2}{3}\right)^{k/2} + 2 \left(\frac{1}{2}\right)^{k/2}
 + \cdots + 2 \left(\frac{1}{129}\right)^{k/2}
 + \frac{2112 \sqrt{33}}{k-3} \left(\frac{8}{33}\right)^{k/2}. \label{r*52bound0}
\end{equation}

We want to show that $|R_{5, 1}^{*}| < 2$ and $|R_{5, 2}^{*}| < 2$. Note that the case for which $(c, d) = \pm (2, 1)$ (resp. $(c, d) = \pm (1, -1)$) yields a bound equal to $2$ for $|R_{5, 1}^{*}|$ (resp. $|R_{5, 2}^{*}|$).

\subsubsection{$\Gamma_0^{*}(7)$}
For $R_{7, 1}^{*}$, we will consider the cases: $N = 2$, $5$, $10$, \ldots, $61$, and $N \geqslant 65$. Then, we have
\begin{equation}
|R_{7, 1}^{*}|
 \leqslant 4
 + 6 \left(\frac{1}{3}\right)^{k/2} + 4 \left(\frac{1}{7}\right)^{k/2}
 + \cdots + 2 \left(\frac{1}{352}\right)^{k/2}
 + \frac{28160}{7(k-3)} \left(\frac{11}{64}\right)^{k/2}. \label{r*71bound0}
\end{equation}

Similarly, for $R_{7, 2}^{*}$, we will consider the cases: $N = 2$, $5$, $10$, \ldots, $89$, and $N \geqslant 97$. For these cases we have
\begin{equation}
|R_{7, 2}^{*}|
 \leqslant 4
 + 2 \left(\frac{1}{2}\right)^{k/2} + 2 \left(\frac{1}{3}\right)^{k/2}
 + \cdots + 2 \left(\frac{1}{571}\right)^{k/2}
 + \frac{62464 \sqrt{6}}{21(k - 3)} \left(\frac{1}{8}\right)^{k/2}. \label{r*72bound0}
\end{equation}

Note that the cases $(c, d) = \pm (2, 1)$ and $\pm (3, 1)$ (resp. $(c, d) = \pm (1, -1)$ and $\pm (3, -1)$) yield a bound equal to $4$ for $|R_{7, 1}^{*}|$ (resp. $|R_{7, 2}^{*}|$).

\subsection{Arguments of some terms}\label{subsec-arg5}

\subsubsection{$\Gamma_0^{*}(5)$}
In the previous subsection, the important point was the fact that the cases $(c, d) = \pm (2, 1)$ and $(c, d) = \pm (1, -1)$ do not yield good bounds for $|R_{5, 1}^{*}|$ and $|R_{5, 2}^{*}|$, respectively.

Let ${\theta_1}' := 2 Arg \left\{ 2 e^{i \theta_1 / 2} + \sqrt{5} e^{- i \theta_1 / 2} \right\}$ and ${\theta_2}' := 2 Arg \left\{ - e^{i \theta_2 / 2} + \sqrt{5} e^{- i \theta_2 / 2} \right\}$, then we have
\begin{equation*}
\tan {\theta_1}' / 2 = - \frac{\sqrt{5} - 2}{\sqrt{5} + 2} \tan \theta_1 / 2, \quad
\tan {\theta_2}' / 2 = - \frac{\sqrt{5} + 1}{\sqrt{5} - 1} \tan \theta_2 / 2.
\end{equation*}

Furthermore, it is easy to show that
\begin{align*}
\tan \left( - \pi + \pi / 2 + \alpha_5 + d_1 (t \pi / k) \right) / 2 &< - \frac{\sqrt{5} - 2}{\sqrt{5} + 2} \tan \left(\pi / 2 + \alpha_5 - (t \pi / k) \right) / 2\\
 &< \tan \left( - \pi + \pi / 2 + \alpha_5 + (t \pi / k) \right) / 2,
\end{align*}
where $d_1 < 1 / (1 + 4 \tan (t / 2) (\pi / k))$. Thus,
\begin{align*}
&\theta_1 = \pi / 2 + \alpha_5 - (t \pi / k)\\
 &\qquad \Rightarrow - \pi + \pi / 2 + \alpha_5 + d_1 (t \pi / k) < {\theta_1}' < - \pi + \pi / 2 + \alpha_5 + (t \pi / k),\\
&k \theta_1 / 2 = k (\pi / 2 + \alpha_5) / 2 - (t / 2) \pi \\
 &\qquad \Rightarrow - (k / 2) \pi + k (\pi / 2 + \alpha_5) / 2 + d_1 (t / 2) \pi < k {\theta_1}' / 2 < - (k / 2) \pi  + k (\pi / 2 + \alpha_5) / 2 + (t / 2) \pi.
\end{align*}

Similarly,
\begin{align*}
&k \theta_2 / 2 = k \alpha_5 / 2 + (t / 2) \pi \\
 &\qquad \Rightarrow - (k / 2) \pi + k \alpha_5 / 2 - (t / 2) \pi < k {\theta_2}' / 2 < - (k / 2) \pi  + k \alpha_5 / 2 - d_2 (t / 2) \pi,
\end{align*}
where $d_2 < 1 / (1 + \tan (t / 2) (\pi / k))$.

Note that
\begin{equation*}
- (k / 2) \pi \equiv
\begin{cases}
0 & (k \equiv 0 \pmod{4})\\
\pi & (k \equiv 2 \pmod{4})
\end{cases} \quad \pmod{2 \pi},
\end{equation*}
and both $d_1$ and $d_2$ tend to $1$ in the limit as $k$ tends to $\infty$, or in the limit as $t$ tends to $0$.

Recall that $\alpha_{5, k} \equiv k (\pi / 2 + \alpha_5) / 2 \pmod{\pi}$, then we can write $k (\pi / 2 + \alpha_5) / 2 = \alpha_{5, k} + m \pi$ for some integer $m$. We define ${\alpha_{5, k}}' \equiv k {\theta_1}' / 2 - m \pi \pmod {2 \pi}$ for $\theta_1 = \pi / 2 + \alpha_5 - (t \pi / k)$.

Similarly, we define $\beta_{5, k} \equiv k \alpha_5 / 2 \pmod{\pi}$ and ${\beta_{5, k}}' \equiv k {\theta_2}' / 2 - (k \alpha_5 / 2 - \beta_{5, k} ) \pmod {2 \pi}$ for $\theta_2 = \alpha_5 + (t \pi / k)$.

\subsubsection{$\Gamma_0^{*}(7)$}\label{subsec-arg7}
Similarly to the previous subsection, we consider the arguments of some terms such that  $(c, d) = \pm (2, 1)$ and $\pm (3, 1)$ for $|R_{7, 1}^{*}|$, and $(c, d) = \pm (1, -1)$ and $\pm (3, -1)$ for $|R_{7, 2}^{*}|$.

Let ${\theta_{1, 1}}' := 2 Arg \left\{ 2 e^{i \theta_1 / 2} + \sqrt{7} e^{- i \theta_1 / 2} \right\}$, ${\theta_{1, 2}}' := 2 Arg \left\{ 3 e^{i \theta_1 / 2} + \sqrt{7} e^{- i \theta_1 / 2} \right\}$,\\
 ${\theta_{2, 1}}' := 2 Arg \left\{ - e^{i \theta_2 / 2} + \sqrt{7} e^{- i \theta_2 / 2} \right\}$, and ${\theta_{2, 2}}' := 2 Arg \left\{ - 3 e^{i \theta_2 / 2} + \sqrt{7} e^{- i \theta_2 / 2} \right\}$.

We have
\begin{align*}
&k \theta_1 / 2 = k (\pi / 2 + \alpha_7) / 2 - (t / 2) \pi \\
 &\qquad \Rightarrow
\begin{cases}
 (k / 3) \pi + k (\pi / 2 + \alpha_7) / 2 + d_{1, 1} (t / 2) \pi\\
 \qquad \qquad < k {\theta_{1, 1}}' / 2 < (k / 3) \pi  + k (\pi / 2 + \alpha_7) / 2 + (3 t / 2) \pi,\\
 - (k / 3) \pi + k (\pi / 2 + \alpha_7) / 2 - t \pi\\
 \qquad \qquad < k {\theta_{1, 2}}' / 2 < - (k / 3) \pi  + k (\pi / 2 + \alpha_7) / 2 - d_{1, 2} (t / 2) \pi,
\end{cases}\\
&k \theta_2 / 2 = k (\alpha_7 - \pi / 6) / 2 + (t / 2) \pi \\
 &\qquad \Rightarrow
\begin{cases}
 - (k / 3) \pi + k (\alpha_7 - \pi / 6) / 2 - (3 t / 4) \pi\\
 \qquad \qquad < k {\theta_{2, 1}}' / 2 < - (k / 3) \pi  + k (\alpha_7 - \pi / 6) / 2 - d_{2, 1} (t / 2) \pi,\\
 (k / 3) \pi + k (\alpha_7 - \pi / 6) / 2 + d_{2, 2} (t / 2) \pi\\
 \qquad \qquad < k {\theta_{2, 2}}' / 2 < (k / 3) \pi  + k (\alpha_7 - \pi / 6) / 2 + (t / 4) \pi,
\end{cases}
\end{align*}
where $d_{1, 1} < 3 / (1 + 2 \sqrt{3} \tan (t / 2) (\pi / k))$, $d_{1, 2} < 2 / (1 + \sqrt{3} \tan (t / 2) (\pi / k))$, $d_{2, 1} < 3 / (2 + \sqrt{3} \tan (t / 2) (\pi / k))$, and $d_{2, 2} < 1 / (2 + 3 \sqrt{3} \tan (t / 2) (\pi / k))$.

Note that
\begin{equation*}
(k / 3) \pi \equiv
\begin{cases}
0 & (k \equiv 0 \pmod{6})\\
2 \pi / 3 & (k \equiv 2 \pmod{6})\\
4 \pi / 3 & (k \equiv 4 \pmod{6})
\end{cases}, \quad
- (k / 3) \pi \equiv
\begin{cases}
0 & (k \equiv 0 \pmod{6})\\
4 \pi / 3 & (k \equiv 2 \pmod{6})\\
2 \pi / 3 & (k \equiv 4 \pmod{6})
\end{cases},
\end{equation*}
modulo $2 \pi$. Furthermore, in the limit as $k$ tends to $\infty$ or in the limit as $t$ tends to $0$,  $d_{1, 1}$ (resp. $d_{1, 2}$, $d_{2, 1}$, and $d_{2, 2}$) tends to $3$ (resp. $2$, $3 / 2$, and $1 / 2$).

Recall that $\alpha_{7, k} \equiv k (\pi / 2 + \alpha_7) / 2 \pmod{\pi}$. Then, we define ${\alpha_{7, k, n}}' \equiv k {\theta_{1, n}}' / 2 - (k (\pi / 2 + \alpha_7) / 2 - \alpha_{7, k} ) \pmod {2 \pi}$ for $n = 1, 2$ and for $\theta_1 = \pi / 2 + \alpha_7 - (t \pi / k)$.

Similarly, we define $\beta_{7, k} \equiv k (\alpha_7 - \pi / 6) / 2 \pmod{\pi}$ and ${\beta_{7, k, n}}' \equiv k {\theta_{2, n}}' / 2 - (k (\alpha_7 - \pi / 6) / 2 - \beta_{7, k} ) \pmod {2 \pi}$ for $n = 1, 2$ and for $\theta_2 = \alpha_7 - \pi / 6 + (t \pi / k)$.

\subsection{Algorithm} \label{subsec-alg}

In this subsection, we consider the bound
\begin{equation}
|R_{p, n}^{*}| < 2 c_0 \quad \text{for every} \; k \geqslant k_0, \label{bound-c0}
\end{equation}
for some $c_0 > 0$ and an even integer $k_0$. Furthermore, we will detail an algorithm that can be used to derive the above bound.

Let $\Lambda$ be an index set, and, applying the RSD method, let us write
\begin{equation*}
|R_{p, n}^{*}| \leqslant 2 {\textstyle \sum_{\lambda \in \Lambda}} \; e_{\lambda}^k \; v_k(c_{\lambda}, d_{\lambda}, \theta),
\end{equation*}
  where the factor ``$2$'' comes from the relation $v_k(c, d, \theta) = v_k(- c, - d, \theta)$. Furthermore, let $I$ be a finite subset of $\Lambda$ such that $e_i^k \; v_k(c_i, d_i, \theta)$ does not tend to $0$ in the limit as $k$ tends to $\infty$ for all $i \in I$, and assume $I \subset \mathbb{N}$. Then, we define $X_i := e_i^{- 2} \; v_k(c_i, d_i, \theta)^{- 2 / k}$ for every $i \in I$.

Assume that for every $i \in I$ and $k \geqslant k_0$ and for some ${c_i}'$ and $u_i$,
\begin{gather*}
\left| Re \left\{ e_i^k \; \left( c_i e^{i \theta / 2} + \sqrt{p} d_i e^{- i \theta / 2} \right)^{- k} \right\} \right| \leqslant {c_i}' X_i^{- k / 2},\\
 X_i^{- k / 2} \geqslant 1 + u_i (\pi / k), \qquad 2 {\textstyle \sum_{\lambda \in \Lambda \setminus I}} \; e_{\lambda}^k \; v_k(c_{\lambda}, d_{\lambda}, \theta) \leqslant b (1 / s)^{k / 2},
\end{gather*}
and let the number $t > 0$ be given.\\

\begin{trivlist}
\item[\bfseries Step 1.]``Determine the number $a_1$.''

First, in order to show the bound (\ref{bound-c0}), we wish to use the bound  
\begin{equation}
{\textstyle \sum_{i \in I}} \; {c_i}' X_i^{- k / 2} < c_0 - a_1 (t \pi / k)^2 \label{bound-a10}
\end{equation}
for every $i \in I$ and $k \geqslant k_0$ and for some $a_1 > 0$.

To show the bound (\ref{bound-c0}) by the above bound (\ref{bound-a10}), we need $b (1 / s)^{k / 2} <  a_1 (t \pi / k)^2$ for every $k \geqslant k_0$. Define $f(k) := s^{k / 2} / b - k^2 / (2 a_1 t^2 \pi^2)$. If we have $k_0 \log s > 4$ and
\begin{equation}
a_1 > (b k_0^2) \; / \; (2 s^{k_0/ 2} t^2 \pi^2), \label{bound-a1}
\end{equation}
then we have $f(k_0) > 0$, $f'(k_0) > 0$, and $f''(k_0) > 0$. In the present paper, we always have $k_0 \log s > 4$. Thus, it is enough to consider the bound (\ref{bound-a1}).\\

\item[\bfseries Step 2.]``Determine the number $c_{0, i}$ and $a_{1, i}$.''

Second, to show bound (\ref{bound-a10}), we wish to use the bounds
\begin{equation}
{c_i}' X_i^{- k / 2} < c_{0, i} - a_{1, i} (t \pi / k)^2 \label{bound-Xi}
\end{equation}
for every $i \in I$ and $k \geqslant k_0$ and for some $c_{0, i} > 0$, $a_{1, i} > 0$.

We determine $c_{0, i}$ and $a_{1, i}$ such that 
\begin{equation}
c_{0, i} > 0, \quad a_{1, i} > 0, \quad c_0 = {\textstyle \sum_{i \in I}} \; c_{0, i} , \quad \text{and} \quad a_1 = {\textstyle \sum_{i \in I}} \; a_{1, i}.
\end{equation}\quad

\item[\bfseries Step 3.]``Determine a discriminant $Y_i$ for every $i \in I$.''

Finally, for the bound (\ref{bound-Xi}), we consider following sufficient conditions:
\begin{equation*}
X_i^{k / 2} > c_i + a_{2, i} (t \pi / k)^2, \quad X_i > a_{3, i} + a_{4, i} (t \pi / k)^2.
\end{equation*}
For the former bound, it is enough to show that
\begin{equation}
c_i = {c_i}' / c_{0, i}, \quad a_{2, i} > c_i^2 (a_{1, i} / {c_i}') \; / \; \{1 - c_i (a_{1, i} / {c_i}') (t \pi / k_0)^2\}, \label{bound-a2}
\end{equation}
while for the latter bound, it is enough to show that
\begin{equation*}
a_{3, i} = c_i^{2 / k_0}, \quad a_{4, i} = \left( (2 a_{2, i}) / (c_i k_0) \right) c_i^{2 / k_0}.
\end{equation*}

Because we have
\begin{equation*}
c_i^{2 / k} \leqslant 1 + 2 (\log c_i) / k + 2 (\log c_i)^2 c_i^{2 / k} / k^2,
\end{equation*}
\begin{align}
X_i - \left( a_{3, i} + a_{4, i} \left( t \frac{\pi}{k} \right)^2 \right)
& \geqslant \frac{1}{k} \left\{ u_i \pi - 2 \log c_i - 2 (\log c_i)^2 c_i^{2 / k} \frac{1}{k_0} - \frac{2 a_{2, i} t^2 \pi^2}{c_i} c_i^{2 / k_0} \frac{1}{k_0^2} \right\} \notag \\
&\quad =: \frac{1}{k} \times Y_i.
\end{align}

In conclusion, if we have $Y_i > 0$, then the bounds (\ref{bound-Xi}), (\ref{bound-a10}), and (\ref{bound-c0}) hold.
\end{trivlist}

Note that the above bounds are sufficient conditions; they are not always necessary.

\section{$\Gamma_0^{*}(5)$ (For Conjecture \ref{conj-g0s5})}
The proof of Conjecture \ref{conj-g0s5} is significantly more difficult than the proof of the theorems for $\Gamma_0^{*}(2)$ and $\Gamma_0^{*}(3)$. The most difficult point concerns the argument $Arg(\rho_{5, 2})$, which is not a rational multiple of $\pi$.

\subsection{All but at most $2$ zeros}\label{subsec-g0s5-ab2}
We have the following lemma:

\begin{lemma}We have the following bounds$:$\\
``We have $|R_{5, 1}^{*}| < 2 \cos({c_0}' \pi)$ for $\theta_1 \in [\pi / 2, \pi / 2 + \alpha_5 - t \pi / k]$''\\
\quad$(1)$ For $k \geqslant 12$, $({c_0}', t) = (1 / 3, 1 / 6)$.\\
\quad$(2)$ For $k \geqslant 58$, $({c_0}', t) = (33 / 80, 9 / 40)$.\\
``We have $|R_{5, 2}^{*}| < 2 \cos({c_0}' \pi)$ for $\theta_2 \in [\alpha_5 + t \pi / k, \pi / 2]$''\\
\quad$(3)$ For $k \geqslant 12$, $({c_0}', t) = (0, 1 / 2)$.\\
\quad$(4)$ For $k \geqslant 22$, $({c_0}', t) = (1 / 3, 1 / 2)$.\\
\quad$(5)$ For $k \geqslant 46$, $({c_0}', t) = (7 / 30, 1 / 5)$.
\label{lem-f5s}
\end{lemma}

\begin{proof}
\noindent
$(3)$ Let $k \geqslant 12$ and $x = \pi / (2 k)$, then $0 \leqslant x \leqslant \pi / 24$, and so $1 - \cos x \geqslant (32/33) x^2$. Thus, we have
\begin{gather*}
\frac{1}{4} |e^{i \theta / 2} - \sqrt{5} e^{-i \theta / 2}|^2
 \geqslant \frac{1}{4} (6 - 2 \sqrt{5} \cos(\alpha_5 + x))
 \geqslant 1 + \frac{16}{11} x^2,\\
\frac{1}{2^k} |e^{i \theta / 2} + \sqrt{5} e^{-i \theta / 2}|^k \geqslant 1 + \frac{96}{11} x^2 \; (k \geqslant 12),\\
2^k \cdot 2 v_k(1, 1, \theta) \leqslant 2 - \frac{288 \pi^2}{\pi^2 + 66} \frac{1}{k^2}.
\end{gather*}
In inequality(\ref{r*52bound0}), replace $2$ with the bound $2 - \frac{288 \pi^2}{\pi^2 + 66} \frac{1}{k^2}$. Then
\begin{equation*}
|R_{5, 2}^{*}|
 \leqslant 2 - \frac{288 \pi^2}{\pi^2 + 66} \frac{1}{k^2}
 +  2 \left(\frac{2}{3}\right)^{k/2}
 + \cdots + 2 \left(\frac{1}{129}\right)^{k/2}
 + \frac{2112 \sqrt{33}}{k-3} \left(\frac{8}{33}\right)^{k/2}.
\end{equation*}
Furthermore,  $(2 / 3)^{k / 2}$ decreases more rapidly in $k$ than $1 / k^2$, and for $k \geqslant 12$, we have
\begin{equation*}
|R_{5, 2}^{*}| \leqslant 1.9821...
\end{equation*}

\noindent
$(1)$, $(2)$, $(4)$, $(5)$ We will use the algorithm in the Subsection \ref{subsec-alg}. Furthermore, we have $X_1 = v_k(2, 1, \theta_1)^{- 2 / k} \geqslant 1 + 4 t (\pi / k)$ in the proof of $(1)$ and $(2)$, and we have $X_1 = (1 / 4) \; v_k(1, -1, \theta_2)^{- 2 / k} \geqslant 1 + t (\pi / k)$ in the proof of $(4)$ and $(5)$. We have $c_0 = c_{0, 1} \leqslant \cos({c_0}' \pi)$. We can show $Y_1 > 0$ for every item.
\end{proof}

When $4 \mid k$, by the valence formula for $\Gamma_0^{*}(5)$ and Proposition \ref{prop-bd_ord_5}, we have at most $k / 4$ zeros on the arc $A_5^{*}$. We have $k / 4 + 1$ {\it integer points} ({\it i.e.} $\cos\left( k \theta / 2 \right) = \pm 1$) in the interval $[\pi / 2, \pi]$. By the above lemma's conditions (1) and (3), we can prove $|R_{5, 1}^{*}| < 2$ or $|R_{5, 2}^{*}| < 2$ at all but at most one integer point. Then, we have all but at most $2$ zeros on $A_5^{*}$.

On the other hand, when $4 \nmid k$, we have at most $(k - 6) / 4$ zeros on the arc $A_5^{*}$. Similarly to the previous case, we have all but at most $2$ zeros on $A_5^{*}$.

Thus, we have the following proposition:
\begin{proposition}
Let $k \geqslant 4$ be an even integer. Then all but at most $2$ of the zeros of $E_{k, 5}^{*}(z)$ in $\mathbb{F}^{*}(5)$ lie on the arc $A_5^{*}$. \label{prop-g0s5-ab2}
\end{proposition}

\subsection{The case $4 \mid k$}\label{subsec-g0s5-40}
For $\pi / 12 < \alpha_{5, k} < 3 \pi / 4$, by Lemma \ref{lem-f5s} (1) and (3), we can prove $|R_{5, 1}^{*}| < 2$ or $|R_{5, 2}^{*}| < 2$ at all of the integer points.

Now, we can write
\begin{gather*}
F_{k, 5, 1}^{*}(\theta_1) = 2 \cos\left( k \theta_1 / 2 \right) + 2 Re(2 e^{- i \theta_1 / 2} + \sqrt{5} e^{i \theta_1 / 2})^{- k} + {R_{5, 1}^{*}}',\\
F_{k, 5, 2}^{*}(\theta_2) = 2 \cos\left( k \theta_2 / 2 \right) + 2^k \cdot 2 Re(e^{- i \theta_2 / 2} - \sqrt{5} e^{i \theta_2 / 2})^{- k} + {R_{5, 2}^{*}}'.
\end{gather*}

For $0 < \alpha_{5, k} < \pi / 12$, the last integer point of $F_{k, 5, 1}^{*}(\theta_1)$ is in the interval $[\pi / 2 + \alpha_5 - \pi / (6 k), \pi / 2 + \alpha_5]$. We have $|{R_{5, 1}^{*}}'| < 2$ for $\theta_1 \in [\pi / 2, \pi / 2 + \alpha_5]$. Furthermore, because $0 < {\alpha_{5, k}}' < \pi / 6$ for $0 < t < 1 / 6$, we have $Sign\{\cos(k \theta_1 / 2)\} = Sign\{Re(2 e^{- i \theta_1 / 2} + \sqrt{5} e^{i \theta_1 / 2})^{- k}\}$ for $\theta_1 \in [\pi / 2 + \alpha_5 - \pi / (6 k), \pi / 2 + \alpha_5]$.

For $3 \pi / 4 < \alpha_{5, k} < \pi$, the first integer point of $F_{k, 5, 2}^{*}(\theta_2)$ is in the interval $[\alpha_5, \alpha_5 + \pi / (2 k)]$. We have $|{R_{5, 2}^{*}}'| < 2$ and $Sign\{\cos(k \theta_2 / 2)\} = Sign\{Re(e^{- i \theta_2 / 2} - \sqrt{5} e^{i \theta_2 / 2})^{- k}\}$ for $\theta_2 \in [\alpha_5, \alpha_5 + \pi / (2 k)]$.

Thus, we have the following proposition:
\begin{proposition}
Let $k \geqslant 4$ be an integer which satisfies $4 \mid k$. Then all of the zeros of $E_{k, 5}^{*}(z)$ in $\mathbb{F}^{*}(5)$ lie on the arc $A_5^{*}$. \label{prop-g0s5-40}
\end{proposition}

\subsection{The case $4 \nmid k$}\label{subsec-g0s5-41}

\subsubsection{\it The case $0 < \alpha_{5, k} < \pi / 2$.}
Now, at most two zeros remain. At the point such that $k \theta_1 / 2 = k (\pi / 2 + \alpha_5) / 2 - \alpha_{5, k} - \pi / 3$, we have $|R_{5, 1}^{*}| < 1$ by Lemma \ref{lem-f5s} (1), and we have $2 \cos(k \theta_1 / 2) = \pm 1$. Then, we have at least one zero between the second to last integer point for $A_{5, 1}^{*}$ and the point $k \theta_1 / 2$. Similarly, by Lemma \ref{lem-f5s} (4), we have at least one zero between the second integer point and the point $k \theta_2 / 2 = k \alpha_5 / 2 + \beta_{5, k} + \pi / 3$.

\subsubsection{\it The case $\pi / 2 < \alpha_{5, k} < \pi$.}
For this case, we expect one more zero between the last integer point for $A_{5, 1}^{*}$ and the first one for $A_{5, 2}^{*}$.  We consider the following cases:

\paragraph{(i)} ``The case $7 \pi / 10 < \alpha_{5, k} < \pi$''\\
\quad$\bullet$ For $3 \pi / 4 < \alpha_{5, k} < \pi$, we can use Lemma \ref{lem-f5s} (1).\\
\quad$\bullet$ For $7 \pi / 10 < \alpha_{5, k} < 3 \pi / 4$, we can use Lemma \ref{lem-f5s} (2).\\
For each case, we consider the point such that $k \theta_1 / 2 = k (\pi / 2 + \alpha_5) / 2 - \alpha_{5, k} + \pi - {c_0}' \pi$. We have $\alpha_{5, k} - \pi + {c_0}' \pi > (t / 2) \pi$ and $|R_{5, 1}^{*}| < 2 \cos({c_0}' \pi)$, and we have $2 \cos(k \theta_1 / 2) = \pm 2 \cos({c_0}' \pi)$. Then, we have at least one zero between the second to  last integer point for $A_{5, 1}^{*}$ and the point $k \theta_1 / 2$.

\paragraph{(ii)} ``The case $\pi / 2 < \alpha_{5, k} < 19 \pi / 30$''\\
\quad$\bullet$ For $\pi / 2 < \alpha_{5, k} < 7 \pi / 12$, we can use Lemma \ref{lem-f5s} (4).\\
\quad$\bullet$ For $7 \pi / 12 < \alpha_{5, k} < 19 \pi / 30$, we can use Lemma \ref{lem-f5s} (5).\\
Similar to the case (i) above, we consider the point such that $k \theta_2 / 2 = k \alpha_5 / 2 - \beta_{5, k} + {c_0}' \pi$ for each case.

\paragraph{(iii)} ``The case $13 \pi / 20 < \alpha_{5, k} < 7 \pi / 10$''

We have $X_1 = v_k(2, 1, \theta_1)^{- 2 / k} \geqslant 1 + 4 t (\pi / k)$, and let $\cos({c_0}' \pi) = - \cos((x / 180) \pi - (t / 2) \pi)$. Then, using the algorithm of Subsection \ref{subsec-alg}, we prove ``For $(x / 180) \pi < \alpha_{5, k} < (y / 180) \pi$, we have $|R_{5, 1}^{*}| < 2 \cos({c_0}' \pi)$ for $\theta_1 = \pi / 2 + \alpha_5 - t \pi / k$.'' for ten cases, namely, $(x, y, t) = (121, 126, 3/20)$, $(120, 121, 1/10)$, $(118.8, 120, 1/10)$, $(118.1, 118.8, 2/25)$, $(117.7, 118.1, 1/15)$, $(117.45, 117.7, 3/50)$, $(117.27, 117.45, 1/20)$, $(117.15, 117.27, 9/200)$, $(117.06, 117.15, 1/25)$, $(117, 117.06, 1/25)$.

For each case, we consider the point such that $k \theta_1 / 2 = k (\pi / 2 + \alpha_5) / 2 - (t / 2) \pi$. We have $\alpha_{5, k} - \pi + {c_0}' \pi > (t / 2) \pi$ and $|R_{5, 1}^{*}| < 2 \cos({c_0}' \pi)$, and we have $|2 \cos(k \theta_1 / 2)| > 2 \cos({c_0}' \pi)$. Then, we have at least one zero between the second to last integer point for $A_{5, 1}^{*}$ and the point $k \theta_1 / 2$.

\paragraph{(iv)} ``The case $19 \pi / 30 < \alpha_{5, k} < 29 \pi / 45$'' 

We have $X_1 = (1 / 4) \; v_k(1, -1, \theta_2)^{- 2 / k} \geqslant 1 + t (\pi / k)$ and $\cos({c_0}' \pi) = \cos((y / 180) \pi - \pi / 2 + (t / 2) \pi)$. Then, we prove ``For $(x / 180) \pi < \alpha_{5, k} < (y / 180) \pi$, we have $|R_{5, 2}^{*}| < 2 \cos({c_0}' \pi)$ for $\theta_2 = \alpha_5 + t \pi / k$.'' for three cases, namely, $(x, y, t) = (114, 115.4, 4/25)$, $(115.4, 115.8, 3/25)$, $(115.8, 116, 1/10)$.

Similar to the case (iii) above, we consider the point such that $k \theta_2 / 2 = k \alpha_5 / 2 + (t / 2) \pi$ for each case.\\

In conclusion, we have the following proposition:
\begin{proposition}
Let $k \geqslant 4$ be an integer which satisfies $4 \nmid k$, and let $\alpha_{5, k} \in [0, \pi]$ be the angle which satisfies $\alpha_{5, k} \equiv k (\pi / 2 + \alpha_5) / 2 \pmod{\pi}$. If we have $\alpha_{5, k} < 29 \pi / 45$ or $13 \pi / 20 < \alpha_{5, k}$, then all of the zeros of $E_{k, 5}^{*}(z)$ in $\mathbb{F}^{*}(5)$ lie on the arc $A_5^{*}$. Otherwise, all but at most one zero of $E_{k, 5}^{*}(z)$ in $\mathbb{F}^{*}(5)$ lie on $A_5^{*}$\label{prop-g0s5-41}
\end{proposition}

\subsection{The remaining case ``$4 \nmid k$ and $29 \pi / 45 < \alpha_{5, k} < 13 \pi / 20$''} \label{subsec-rest5}
In the previous subsection, we left one zero between the last integer point for $A_{5, 1}^{*}$ and the first one for $A_{5, 2}^{*}$ for the case of ``$4 \nmid k$ and $29 \pi / 45 < \alpha_{5, k} < 13 \pi / 20$''. For the cases of $13 \pi / 20 < \alpha_{5, k} < 7 \pi / 10$ and $19 \pi / 30 < \alpha_{5, k} < 29 \pi / 45$, the width $|x - y|$ becomes smaller as the intervals of the bounds approach the interval $[29 \pi / 45, 13 \pi / 20]$. It seems that the width $|x - y|$ needs to be smaller still if we are to prove Conjecture \ref{conj-g0s5} for the remaining interval $[29 \pi / 45, 13 \pi / 20]$. Furthermore, we may need to split infinite cases $(x, y)$ such as we saw in the previous subsection. Thus, we cannot prove the conjecture for this remaining case in a similar manner. However, when $k$ is large enough, there is a possibility that we can prove the conjecture for this remaining case.

Let $29 \pi / 45 < \alpha_{5, k} < 13 \pi / 20$, and let $t > 0$ be small enough. Then, we have $\pi / 2 < \alpha_{5, k} - (t / 2) \pi < \pi$ and $3 \pi / 2 < \pi + \alpha_{5, k} + d_1 (t / 2) \pi < {\alpha_{5, k}}' < \pi + \alpha_{5, k} + (t / 2) \pi < 2 \pi$. Moreover, we can easily show that $1 + 4 t (\pi / k) \leqslant v_k(2, 1, \theta_1)^{- 2 / k} \leqslant e^{4 t (\pi / k)}$. Thus, we have
\begin{align*}
&- \cos(\alpha_{5, k} - (t / 2) \pi) - \cos(\pi + \alpha_{5, k} + d_1 (t / 2) \pi) \cdot e^{- 2 \pi t}\\
 &\qquad > |\cos(k \theta_1 / 2)| - \left| Re\left\{ \left( 2 e^{i \theta_1 / 2} + \sqrt{5} e^{- i \theta_1 / 2} \right)^{- k} \right\} \right|\\
 &\qquad > - \cos(\alpha_{5, k} - (t / 2) \pi) - \cos(\pi + \alpha_{5, k} + (t / 2) \pi) \cdot (1 + 4 t (\pi / k))^{- k / 2}.
\end{align*}
We denote the upper bound by $A$ and the lower bound by $B$. Furthermore, we define $A' := A / \cos(\pi + \alpha_{5, k} + d_1 (t / 2) \pi)$ and $B' := B / \cos(\pi + \alpha_{5, k} + (t / 2) \pi)$. First, we have $A |_{t = 0} = B |_{t = 0} = 0$. Second, we have $\frac{\partial}{\partial t} A' |_{t = 0} = \frac{\partial}{\partial t} B' |_{t = 0} = \pi (\tan \alpha_{k, 5} + 2)$. Finally, we have $B > 0$ if $\alpha_{5, k} > \pi - \alpha_5$, and we have $A < 0$ if $\alpha_{5, k} < \pi - \alpha_5$ for small enough $t$.

Similarly, we consider the lower and the upper bounds of $|\cos(k \theta_2 / 2)| - | Re\{ 2^k ( e^{i \theta_2 / 2} - \sqrt{5} e^{- i \theta_2 / 2} )^{- k} \} |$.  The lower bound is positive if $\alpha_{5, k} < \pi - \alpha_5$, while the upper bound is negative if $\alpha_{5, k} > \pi - \alpha_5$ for small enough $t$.

In conclusion, if $4 \nmid k$ is large enough, then $|{R_{5, 1}^{*}}'|$ and $|{R_{5, 2}^{*}}'|$ are small enough, and then we have one more zero on the arc $A_{5, 1}^{*}$ when $\alpha_{5, k} > \pi - \alpha_5$, and one more zero on the arc $A_{5, 2}^{*}$ when $\alpha_{5, k} < \pi - \alpha_5$. However, if $k$ is small, a method of proving the conjecture for this case is not clear.

\section{$\Gamma_0^{*}(7)$ (For Conjecture \ref{conj-g0s7})}
Similar to the case of $\Gamma_0^{*}(5)$, to prove Conjecture \ref{conj-g0s7} is also difficult. The most difficult point is again the argument $Arg(\rho_{7, 2})$.

\subsection{All but at most $2$ zeros}\label{subsec-g0s7-ab2}
We have the following lemma:

\begin{lemma}We have the following bounds$:$\\
``We have $|R_{7, 1}^{*}| < 2 \cos({c_0}' \pi)$ for $\theta_1 \in [\pi / 2, \pi / 2 + \alpha_7 - t \pi / k]$''\\
\quad$(1)$ For $k \geqslant 10$, $({c_0}', t) = (1 / 3, 1 / 3)$.\\
\quad$(2)$ For $k \geqslant 80$, $({c_0}', t) = (41 / 100, 8 / 25)$.\\
\quad$(3)$ For $k \geqslant 22$, $({c_0}', t) = (13 / 36, 1 / 3)$.\\
``We have $|R_{7, 2}^{*}| < 2 \cos({c_0}' \pi)$ for $\theta_2 \in [\alpha_7 - \pi / 6 + t \pi / k, \pi / 2]$''\\
\quad$(4)$ For $k \geqslant 8$, $({c_0}', t) = (1 / 6, 1 / 2)$.\\
\quad$(5)$ For $k = 26$, $k \geqslant 44$, $({c_0}', t) = (1 / 3, 2 / 3)$.\\
\quad$(6)$ For $k \geqslant 70$, $({c_0}', t) = (1 / 4, 1 / 2)$.\\
\quad$(7)$ For $k \geqslant 200$, $({c_0}', t) = (5 / 18, 1 / 2)$.\label{lem-f7s}
\end{lemma}

\begin{proof}
We will use the algorithm given in Subsection \ref{subsec-alg}. Furthermore, we have $X_1 = v_k(2, 1, \theta_1)^{- 2 / k} \geqslant 1 + 2 \sqrt{3} t (\pi / k)$ and $X_2 = v_k(3, 1, \theta_1)^{- 2 / k} \geqslant 1 + 3 \sqrt{3} t (\pi / k)$ in the proofs of (1), (2), and (3), and we have $X_1 = (1 / 4) \; v_k(1, -1, \theta_2)^{- 2 / k} \geqslant 1 + (\sqrt{3} / 2) t (\pi / k)$ and $X_2 = (1 / 4) \; v_k(3, -1, \theta_2)^{- 2 / k} \geqslant 1 + (3 \sqrt{3} / 2) t (\pi / k)$ in the proofs of (4), (5), (6), and (7). We also have $c_0 \leqslant \cos({c_0}' \pi)$. We can show $X_i - (a_{3, i} + a_{4, i} (t \pi / k)^2) > 0$ in the algorithm given in Subsection \ref{subsec-alg} for the case of ``(4), $k = 8$'' and ``(5), $k = 26$''. For the other cases, we can show $Y_1 > 0$ and $Y_2 > 0$.
\end{proof}

Similarly to Proposition \ref{prop-g0s5-ab2}, by the above lemma's conditions (1), (4) and (5), we have the following proposition:
\begin{proposition}
Let $k \geqslant 4$ be an even integer. Then all but at most $2$ of the zeros of $E_{k, 7}^{*}(z)$ in $\mathbb{F}^{*}(7)$ lie on the arc $A_7^{*}$. \label{prop-g0s7-ab2}
\end{proposition}

\subsection{The case $6 \mid k$}\label{subsec-g0s7-60}
We can write
\begin{align*}
F_{k, 7, 1}^{*}(\theta_1) = &2 \cos\left( k \theta_1 / 2 \right) + 2 Re(2 e^{- i \theta_1 / 2} + \sqrt{7} e^{i \theta_1 / 2})^{- k}\\
 &\qquad+ 2 Re(3 e^{- i \theta_1 / 2} + \sqrt{7} e^{i \theta_1 / 2})^{- k} + {R_{7, 1}^{*}}',\\
F_{k, 7, 2}^{*}(\theta_2) = &2 \cos\left( k \theta_2 / 2 \right) + 2^{k + 1} Re(e^{- i \theta_2 / 2} - \sqrt{7} e^{i \theta_2 / 2})^{- k}\\
 &\qquad+ 2^{k + 1} Re(3 e^{- i \theta_2 / 2} - \sqrt{7} e^{i \theta_2 / 2})^{- k} + {R_{7, 2}^{*}}'.
\end{align*}
Similarly to Subsection \ref{prop-g0s5-40}, we consider the signs of some of these terms.

For $0 < \alpha_{7, k} < \pi / 8$, we have $|{R_{7, 1}^{*}}'| < 2$ and $Sign\{\cos(k \theta_1 / 2)\} = Sign\{Re(2 e^{- i \theta_1 / 2} + \sqrt{7} e^{i \theta_1 / 2})^{- k}\} = Sign\{Re(3 e^{- i \theta_1 / 2} + \sqrt{7} e^{i \theta_1 / 2})^{- k}\}$ for $\theta_1 \in [\pi / 2 + \alpha_7 - \pi / (8 k), \pi / 2 + \alpha_7]$.

For $\pi / 8 < \alpha_{7, k} < \pi / 6$ or $\pi / 4 < \alpha_{7, k} < 5 \pi / 6$, we can use Lemma \ref{lem-f7s} (1) and (4).

For $5 \pi / 6 < \alpha_{7, k} < \pi$, we have$|{R_{7, 2}^{*}}'| < 2$ and $Sign\{\cos(k \theta_2 / 2)\} = Sign\{Re(e^{- i \theta_2 / 2} - \sqrt{7} e^{i \theta_2 / 2})^{- k}\} = Sign\{Re(3 e^{- i \theta_2 / 2} - \sqrt{7} e^{i \theta_2 / 2})^{- k}\}$ for $\theta_2 \in [\alpha_7 - \pi / 6 + \pi / (6 k), \pi / 2]$.

Thus, we have the following proposition:
\begin{proposition}
Let $k \geqslant 4$ be an integer which satisfies $6 \mid k$. Then all of the zeros of $E_{k, 7}^{*}(z)$ in $\mathbb{F}^{*}(7)$ lie on the arc $A_7^{*}$. \label{prop-g0s7-60}
\end{proposition}

\subsection{The case $k \equiv 2 \pmod{6}$}\label{subsec-g0s7-61}
We can prove this case in a similar way to that of Subsection \ref{subsec-g0s5-41}. 

\subsubsection{\it The case $0 < \alpha_{7, k} < 2 \pi / 3$.}
We can use Lemma \ref{lem-f7s} (1), (4), and (5). When $\alpha_{7, k} < \pi / 6$, we consider the point $k \theta_1 / 2 =  k (\pi / 2 + \alpha_7) / 2 - \alpha_{7, k} - \pi / 3$ instead of the last integer point for $A_{7, 1}^{*}$. Similarly, instead of the first integer point for $A_{7, 2}^{*}$, we consider the points $k \theta_2 / 2 =  k (\alpha_7 - \pi / 6) / 2 + (\pi - \beta_{7, k}) + \pi / 6$ and $k (\alpha_7 - \pi / 6) / 2 + (\pi - \beta_{7, k}) + \pi / 3$ for $5 \pi / 12 < \alpha_{7, k} < 7 \pi / 12$ and $7 \pi / 12 < \alpha_{7, k} < 2 \pi / 3$, respectively.

\subsubsection{\it The case $2 \pi / 3 < \alpha_{7, k} < \pi$.}
For this case, we expect one more zero between the last integer point for $A_{7, 1}^{*}$ and the first one for $A_{7, 2}^{*}$. Then, we consider the following cases:

\paragraph{(i)} ``The case $3 \pi / 4 < \alpha_{7, k} < \pi$''

For $3 \pi / 4 < \alpha_{7, k} < 5 \pi / 6$ and $5 \pi / 6 < \alpha_{7, k} < \pi$, we can use Lemma \ref{lem-f7s} (2) and (1), respectively.

\paragraph{(ii)} ``The case $3217 \pi / 4500 < \alpha_{7, k} < 3 \pi / 4$''

We have $X_1 = v_k(2, 1, \theta_1)^{- 2 / k} \geqslant 1 + 2 \sqrt{3} t (\pi / k)$ and $X_2 = v_k(3, 1, \theta_1)^{- 2 / k} \geqslant 1 + 3 \sqrt{3} t (\pi / k)$, and let $\cos({c_0}' \pi) = - \cos((x / 180) \pi - (t / 2) \pi)$. Then, by the algorithm given in Subsection \ref{subsec-alg}, we prove ``For $(x / 180) \pi < \alpha_{7, k} < (y / 180) \pi$, we have $|R_{7, 1}^{*}| < 2 \cos({c_0}' \pi)$ for $\theta_1 = \pi / 2 + \alpha_7 - t \pi / k$.'' for nine cases, namely, $(x, y, t) = (131.5, 135, 1/4)$, $(130.1, 131.5, 83/400)$, $(129.5, 130.1, 7/40)$, $(129.18, 129.5, 47/300)$, $(129, 129.18, 71/500)$, $(128.86, 129, 263/2000)$, $(128.77, 128.86, 61/500)$, $(128.71, 128.77, 143/1250)$,\\ $(128.68, 128.71, 109/1000)$.

\paragraph{(iii)} ``The case $2 \pi / 3 < \alpha_{7, k} < 266 \pi / 375$'' 

We have $X_1 = (1 / 4) \; v_k(1, -1, \theta_2)^{- 2 / k} \geqslant 1 + (\sqrt{3} / 2) t (\pi / k)$ and $X_2 = (1 / 4) \; v_k(3, -1, \theta_2)^{- 2 / k} \geqslant 1 + (3 \sqrt{3} / 2) t (\pi / k)$, and let $\cos({c_0}' \pi) = \cos((y / 180) \pi - 2 \pi / 3 - (t / 2) \pi)$. Then, we prove ``For $(x / 180) \pi < \alpha_{7, k} < (y / 180) \pi$, we have $|R_{7, 2}^{*}| < 2 \cos({c_0}' \pi)$ for $\theta_2 = \alpha_7 - \pi / 6 + t \pi / k$.'' for four cases, namely, $(x, y, t) = (120, 126.7, 93/200)$, $(126.7, 127.3, 17/50)$, $(127.3, 127.63, 22/75)$, $(127.63, 127.68, 13/50)$.\\

In conclusion, we have the following proposition:
\begin{proposition}
Let $k \geqslant 4$ be an integer which satisfies $k \equiv 2 \pmod{6}$, and let $\alpha_{7, k} \in [0, \pi]$ be the angle which satisfies $\alpha_{7, k} \equiv k (\pi / 2 + \alpha_7) / 2 \pmod{\pi}$. If we have $\alpha_{7, k} < 266 \pi / 375$ or $3217 \pi / 4500 < \alpha_{7, k}$, then all of the zeros of $E_{k, 7}^{*}(z)$ in $\mathbb{F}^{*}(7)$ lie on the arc $A_7^{*}$. Otherwise, all but at most one zero of $E_{k, 7}^{*}(z)$ in $\mathbb{F}^{*}(7)$ lie on $A_7^{*}$\label{prop-g0s7-61}
\end{proposition}

\subsection{The case $k \equiv 4 \pmod{6}$}\label{subsec-g0s7-62}
With the exception of some specific cases, we can prove this case in a similar way to the proof of Subsection \ref{subsec-g0s5-41} and the previous subsection.

\subsubsection{\it The case $0 < \alpha_{7, k} < \pi / 3$.}
We can use Lemma \ref{lem-f7s} (1) for the case $\alpha_{7, k} < \pi / 6$, and we can use Lemma \ref{lem-f7s} (4) and (6) for the cases $0 < \alpha_{7, k} < \pi / 4$ and $\pi / 4 < \alpha_{7, k} < \pi / 3$, respectively.

\subsubsection{\it The case $\pi / 3 < \alpha_{7, k} < \pi$.}

\paragraph{(i)} ``The case $3 \pi / 4 < \alpha_{7, k} < \pi$''

For $5 \pi / 6 < \alpha_{7, k} < \pi$ (resp. $29 \pi / 36 < \alpha_{7, k} < 5 \pi / 6$, $3 \pi / 4 < \alpha_{7, k} < 5 \pi / 6$), we can use Lemma \ref{lem-f7s} (1) (resp. (3), (2)).

\paragraph{(ii)} ``The case $\pi / 3 < \alpha_{7, k} < 13 \pi / 36$'' We can use Lemma \ref{lem-f7s} (7).

\paragraph{(iii)} ``The case $2 \pi / 3 < \alpha_{7, k} < 3 \pi / 4$''

We define $\cos({c_0}' \pi) = - \cos((x / 180) \pi - (t / 2) \pi)$. Then, we prove ``For $(x / 180) \pi < \alpha_{7, k} < (y / 180) \pi$, we have $|R_{7, 1}^{*}| < 2 \cos({c_0}' \pi)$ for $\theta_1 = \pi / 2 + \alpha_7 - t \pi / k$.'' for two cases, namely, $(x, y, t) = (127.6, 135, 59 / 250)$, $(120, 127.6, 1 / 4)$.

\paragraph{(iv)} ``The case $13 \pi / 36 < \alpha_{7, k} < 5 \pi / 9$'' 

We define $\cos({c_0}' \pi) = \cos((y / 180) \pi - \pi / 3 + (t / 2) \pi)$. Then, we prove ``For $(x / 180) \pi < \alpha_{7, k} < (y / 180) \pi$, we have $|R_{7, 2}^{*}| < 2 \cos({c_0}' \pi)$ for $\theta_2 = \alpha_7 - \pi / 6 + t \pi / k$.'' for two cases, namely, $(x, y, t) = (65, 90, 2 / 5)$, $(90, 100, 2 / 5)$.\\

Now, we can write
\begin{gather*}
F_{k, 7, 1}^{*}(\theta_1) = 2 \cos\left( k \theta_1 / 2 \right) + 2 Re(3 e^{- i \theta_1 / 2} + \sqrt{7} e^{i \theta_1 / 2})^{- k} + {R_{7, 1}^{*}}'',\\
F_{k, 7, 2}^{*}(\theta_2) = 2 \cos\left( k \theta_2 / 2 \right) + 2^k \cdot 2 Re(3 e^{- i \theta_2 / 2} - \sqrt{7} e^{i \theta_2 / 2})^{- k} + {R_{7, 2}^{*}}''.
\end{gather*}

\paragraph{(v)} ``The case $73 \pi / 120 < \alpha_{7, k} < 2 \pi / 3$'' 

We have $X_1 = v_k(2, 1, \theta_1)^{- 2 / k} \geqslant 1 + 2 \sqrt{3} t (\pi / k)$ and $X_2 = v_k(3, 1, \theta_1)^{- 2 / k} \leqslant e^{- 3 \sqrt{3} (t / 2) \pi}$. Then,   $Sign\{\cos\left( k \theta_1 / 2 \right)\} = Sign\{Re(3 e^{- i \theta_1 / 2} + \sqrt{7} e^{i \theta_1 / 2})^{- k}\}$, and we prove ``For $(x / 180) \pi < \alpha_{7, k} < (y / 180) \pi$, we have $|{R_{7, 1}^{*}}''| < |2 \cos\left( k \theta_1 / 2 \right) + 2 Re(3 e^{- i \theta_1 / 2} + \sqrt{7} e^{i \theta_1 / 2})^{- k}|$ for $\theta_1 = \pi / 2 + \alpha_7 - t \pi / k$.'' for four cases, namely, $(x, y) = (111.6, 120, 23 / 150)$, $(110.1, 111.6, 1 / 10)$, $(109.65, 110.1, 43 / 625)$, $(109.5, 109.65, 21 / 400)$.

For each case, we consider the point such that $k \theta_1 / 2 = k (\pi / 2 + \alpha_7) / 2 - (t / 2) \pi$.  We can show $Sign\{\cos\left( k \theta_1 / 2 \right)\} = Sign\{F_{k, 7, 1}^{*}(\theta_1)\}$ , and then we have at least one zero between the second to last integer point for $A_{7, 1}^{*}$ and the point $k \theta_1 / 2$.

\paragraph{(vi)} ``The case $5 \pi / 9 < \alpha_{7, k} < 217 \pi / 360$'' 

We have $X_1 = (1 / 4) \; v_k(1, -1, \theta_2)^{- 2 / k} \geqslant 1 + t (\pi / k)$ and $X_2 = (1 / 4) \; v_k(3, -1, \theta_2)^{- 2 / k} \leqslant e^{- (3 \sqrt{3} / 2) (t / 2) \pi}$. Then, we have $Sign\{\cos\left( k \theta_2 / 2 \right)\} = Sign\{Re(3 e^{- i \theta_2 / 2} - \sqrt{7} e^{i \theta_2 / 2})^{- k}\}$, and we prove ``For $(x / 180) \pi < \alpha_{7, k} < (y / 180) \pi$, we have $|{R_{7, 2}^{*}}''| < |2 \cos\left( k \theta_2 / 2 \right) + 2^k \cdot 2 Re(3 e^{- i \theta_2 / 2} - \sqrt{7} e^{i \theta_2 / 2})^{- k}|$ for $\theta_2 = \alpha_7 - \pi / 6 + t \pi / k$.'' for five cases, namely, $(x, y) = (100, 106, 3/10)$, $(106, 107.7, 11 / 50)$, $(107.7, 108.21, 33 / 200)$, $(108.21, 108.42, 2 / 15)$, $(108.42, 108.5, 113 / 1000)$.\\

In conclusion, we have the following proposition:
\begin{proposition}
Let $k \geqslant 4$ be an integer which satisfies $k \equiv 4 \pmod{6}$, and let $\alpha_{7, k} \in [0, \pi]$ be the angle which satisfies $\alpha_{7, k} \equiv k (\pi / 2 + \alpha_7) / 2 \pmod{\pi}$. If we have $\alpha_{7, k} < 217 \pi / 360$ or $73 \pi / 120 < \alpha_{7, k}$, then all of the zeros of $E_{k, 7}^{*}(z)$ in $\mathbb{F}^{*}(7)$ lie on the arc $A_7^{*}$. Otherwise, all but at most one zero of $E_{k, 7}^{*}(z)$ in $\mathbb{F}^{*}(7)$ are on $A_7^{*}$\label{prop-g0s7-62}
\end{proposition}

\subsection{The remaining cases ``$k \equiv 2 \pmod{6}$, $266 \pi / 375 < \alpha_{7, k} < 3217 \pi / 4500$'' and ``$k \equiv 4 \pmod{6}$, $217 \pi / 360 < \alpha_{7, k} < 73 \pi / 120$'' } \label{subsec-rest7}
Similar the problem described in Subsection \ref{subsec-rest5}, it is difficult to prove Conjecture \ref{conj-g0s7} for the remaining cases. However, when $k$ is large enough, we have the following observation.

\subsubsection{\bfseries The case ``$k \equiv 2 \pmod{6}$ and $266 \pi / 375 < \alpha_{7, k} < 3217 \pi / 4500$''}
Let $t > 0$ be small enough, then we have $\pi / 2 < \alpha_{7, k} - (t / 2) \pi < \pi$, $\pi < 2 \pi / 3 + \alpha_{7, k} + d_{1, 1} (t / 2) \pi < {\alpha_{7, k, 1}}' < 2 \pi / 3 + \alpha_{7, k} + (t / 2) \pi < 3 \pi / 2$, and $2 \pi < 4 \pi / 3 + \alpha_{7, k} - t \pi < {\alpha_{7, k, 2}}' < 4 \pi / 3 + \alpha_{7, k} - d_{1, 2} (t / 2) \pi < 5 \pi / 2$. Thus, we have
\begin{allowdisplaybreaks}
\begin{align*}
&- \cos(\alpha_{7, k} - (t / 2) \pi)
 - \cos(2 \pi / 3 + \alpha_{7, k} + d_{1, 1} (t / 2) \pi) \cdot (1 + 2 \sqrt{3} t (\pi / k))^{- k / 2}\\
  &\qquad \qquad \qquad - \cos(4 \pi / 3 + \alpha_{7, k} - d_{1, 2} (t / 2) \pi) \cdot e^{- (3 \sqrt{3} / 2) \pi t}\\
 &\qquad > |\cos(k \theta_1 / 2)|
  + \left| Re\left\{ \left( 2 e^{i \theta_1 / 2} + \sqrt{7} e^{- i \theta_1 / 2} \right)^{- k} \right\} \right|\\
  &\qquad \qquad \qquad \qquad - \left| Re\left\{ \left( 3 e^{i \theta_1 / 2} + \sqrt{7} e^{- i \theta_1 / 2} \right)^{- k} \right\} \right|\\
 &\qquad > - \cos(\alpha_{7, k} - (t / 2) \pi)
 - \cos(2 \pi / 3 + \alpha_{7, k} + (3 t / 2) \pi) \cdot e^{- \sqrt{3} \pi t}\\
  &\qquad \qquad \qquad \qquad - \cos(4 \pi / 3 + \alpha_{7, k} - t \pi) \cdot (1 + 3 \sqrt{3} t (\pi / k))^{- k / 2}.
\end{align*}
\end{allowdisplaybreaks}
We denote the upper bound by $A$ and the lower bound by $B$. First, we have $A |_{t = 0} = B |_{t = 0} = 0$ and $\frac{\partial}{\partial t} A |_{t = 0} = \frac{\partial}{\partial t} B |_{t = 0} = 0$. Second, let $C = (5 \sqrt{3} / 2) \pi^2 (- \cos \alpha_{7, k}) (\tan \alpha_{7, k} + 11 / (5 \sqrt{3}))$, then we have $\frac{\partial^2}{\partial t^2} A |_{t = 0} = C + 6 \pi^2 (- \cos(2 \pi / 3 + \alpha_{7, k})) / k$ and $\frac{\partial^2}{\partial t^2} B |_{t = 0} = C - (27 / 2) \pi^2 \cos(4 \pi / 3 + \alpha_{7, k}) / k$. Finally, we have $B > 0$ if $\alpha_{7, k} > 3 \pi / 2 - 2 \alpha_7$, and we have $A < 0$ if $\alpha_{7, k} < 3 \pi / 2 - 2 \alpha_7$ for large enough $k$ and small enough $t$.

Similarly, we consider the lower and the upper bounds of $|\cos(k \theta_2 / 2)|  - | Re\{ 2^k \cdot ( e^{i \theta_2 / 2} - \sqrt{7} e^{- i \theta_2 / 2} )^{- k} \} |   - | Re\{ 2^k \cdot ( 3 e^{i \theta_2 / 2} - \sqrt{7} e^{- i \theta_2 / 2} )^{- k} \} |$. The lower bound is positive if $\alpha_{7, k} < 3 \pi / 2 - 2 \alpha_7$, while the upper bound is negative if $\alpha_{7, k} > 3 \pi / 2 - 2 \alpha_7$ for large enough $k$ and for small enough $t$.

In conclusion, if $k$ is large enough, then $|{R_{7, 1}^{*}}'|$ and $|{R_{7, 2}^{*}}'|$ is small enough, and then we have one more zero on the arc $A_{7, 1}^{*}$ when $\alpha_{7, k} > 3 \pi / 2 - 2 \alpha_7$, and we have one more zero on the arc $A_{7, 2}^{*}$ when $\alpha_{7, k} < 3 \pi / 2 - 2 \alpha_7$. However, if $k$ is small, a method of proving the conjecture for this case is not clear.

\subsubsection{\bfseries The case ``$k \equiv 4 \pmod{6}$ and $217 \pi / 360 < \alpha_{7, k} < 73 \pi / 120$''}
Let $t > 0$ be small enough. Similar to Subsection \ref{subsec-rest5}, if $k$ is large enough, then we have one more zero on the arc $A_{7, 1}^{*}$ when $\alpha_{7, k} > \pi - \alpha_7$, and we have one more zero on the arc $A_{7, 2}^{*}$ when $\alpha_{7, k} < \pi - \alpha_7$.\\

\noindent
{\itshape Acknowledgement.} I would like to thank Professor Eiichi Bannai for suggesting these problems.

\quad\\

\begin{flushright}\makebox{
\begin{minipage}{2.8in}{\it
\begin{center}
Junichi Shigezumi\\
Graduate School of Mathematics\\
Kyushu University\\
Hakozaki 6-10-1 Higashi-ku\\
Fukuoka 812-8581\\
Japan\\
$($E-mail: j.shigezumi@math.kyushu-u.ac.jp$)$
\end{center}
}\end{minipage}}\end{flushright}

\begin{thebibliography}{RSD}

\bibitem[G]{G}
J. Getz, {\it A generalization of a theorem of Rankin and Swinnerton-Dyer on zeros of modular forms.}, Proc. Amer. Math. Soc., {\bfseries 132}(2004), No. 8, 2221-2231.

\bibitem[K]{K}
A. Krieg, {\it Modular Forms on the Fricke Group.}, Abh. Math. Sem. Univ. Hamburg, {\bfseries 65}(1995), 293-299.

\bibitem[MNS]{MNS}
T. Miezaki, H. Nozaki, J. Shigezumi, {\it On the zeros of Eisenstein series for $\Gamma_0^* (2)$ and $\Gamma_0^* (3)$}, to appear.

\bibitem[Q]{Q}
H. -G. Quebbemann, {\it Atkin-Lehner eigenforms and strongly modular lattices}, Enseign. Math. (2), {\bfseries 4}3(1997), No. 1-2, 55-65.

\bibitem[RSD]{RSD}
F. K. C. Rankin, H. P. F. Swinnerton-Dyer, {\it On the zeros of Eisenstein Series}, Bull. London Math. Soc., {\bfseries 2}(1970), 169-170.

\bibitem[SE]{SE}
J. -P. Serre, {\it A Corse in Arithmetic}, Graduate Texts in Mathematics, No. 7, Springer-Verlag, New York-Heidelberg, 1973. (Translation of {\it Cours d'arithm\'etique $($French$)$}, Presses Univ. France, Paris, 1970.)

\bibitem[SG]{SG}
G. Shimura, {\it Introduction to the arithmetic theory of automorphic functions},  Kan\^{o} Memorial Lectures, No. 1. Publ. Math. Soc. Japan, No. 11. Iwanami Shoten Publishers, Tokyo; Princeton Univ. Press, Princeton, 1971.

\bibitem[SH]{SH}
H. Shimizu, {\it Hokei kansu. I-III. $($Japanese$)$ $[$Automorphic functions. I-III$]$}, Iwanami Shoten Kiso Sugaku [Iwanami Lectures on Fundamental Mathematics] 8, Iwanami Shoten Publishers, Tokyo, 1977-1978.

\bibitem[SJ]{SJ}
J. Shigezumi, {\it A detailed note on the zeros of Eisenstein series for $\Gamma_0^{*}(5)$ and $\Gamma_0^{*}(7)$}, arXiv:math.NT/0607247.

\end{thebibliography}
\end{document}